\documentclass[12pt,a4paper]{article}
\usepackage[a4paper,top=2.5cm,bottom=2.5cm,left=2.5cm,right=2.5cm]{geometry}
\usepackage{amsmath}
 \usepackage{graphicx}
\usepackage{amsthm}
\usepackage{color}
\usepackage{fancyhdr}
\usepackage{pdfpages}
\usepackage{amssymb}
\usepackage{paralist}
\usepackage{extarrows}
\usepackage{setspace}
\usepackage{braket} 
\usepackage[colorlinks=true,citecolor=blue,linkcolor=blue,anchorcolor=blue,urlcolor=blue]{hyperref}
\allowdisplaybreaks
\numberwithin{equation}{section}
\newtheorem{theorem}{Theorem}[section]
\newtheorem{lemma}{Lemma} [section]

\newtheorem{definition}[theorem]{Definition}
\newtheorem{remark}{Remark}[section]
\begin{document}
\title{Global existence and finite time blow-up  for the heat flow of  H-system with
 constant mean curvature}
\author{Fei Fang \footnote{E-mail: fangfei68@163.com} and Yannan Liu \\  \footnotesize  \emph{School of Mathematics and Statistics,%
Beijing Technology and Business University, Beijing 100048, China}}
\maketitle
\noindent \textbf{\textbf{Abstract:}}  In this paper, we use the modified potential
 well method to  study the long time behaviors of solutions to the heat flow of
  H-system in a bounded smooth domain of $R^2$. Global existence and finite time
   blowup of solutions are proved when the initial energy is in three cases.
   When the initial energy is low or critical, we not only give a threshold
   result for  the global existence and blowup of solutions, but also
   obtain the decay rate of the $L^2$ norm for global solutions. When
   the initial energy is high, sufficient conditions for the global
   existence and blowup of solutions are also provided. We extend
   the recent results which were obtained in \cite{r4}.


\section{Introduction}
Let $\Omega \subset \mathbb{R}^{2}$ be a bounded domain with smooth boundary $\partial \Omega$, and let $H$ be a
bounded Lipschitz function on $\mathbb{R}^{3} .$ A map $u \in C^{2}\left(\Omega, \mathbb{R}^{3}\right)$ is called an $H$-surface
(parametrized over $\Omega$ ) if $u$ satisfies
\begin{equation}\label{e0}
  \Delta u=2 H(u) u_{x} \wedge u_{y},
\end{equation}
where $\wedge$ denotes the wedge product of $\mathbb{R}^{3}.$
System of the general form (\ref{e0}) arises from  differential geometry and in the calculus of variation.
If  $u$ is a conformal representation
of a surface $S$ in $\mathbb{R}^{3},$ i.e., $u_{x} \cdot u_{y}=0=\left|u_{x}\right|^{2}-\left|u_{y}\right|^{2},$ then $H(u)$ is the mean
curvature of $S$ at the point $u.$  For $H(u)\equiv const$, the weak solutions of the Dirichlet problem associated to (\ref{e0}) correspond to critical
points of the energy functional,
$$e(u)=\frac{1}{3}\int_{\Omega} |\nabla u|^2,$$
under the constraint that the volume functional
$$V_H(u)=\frac{2}{3} \int_{\Omega} H(u) u \cdot u_{x} \wedge u_{y},$$
is a given constant.

Starting with the pioneering works  of Wente \cite{r7} in $1969,$
a very large
amount of literature has been devoted to   system \ref{e0}.
 System \ref{e0}
with constant mean curvature $H$ has been extensively studied by Wente \cite{r7}, Hildebrandt \cite{r8},
Struwe \cite{r9}, and Brezis and Coron \cite{r10,r10}.
 Hildebrant \cite{r8} considered the Plateau problem for
surfaces of constant mean curvature, Brezis and Coron \cite{r10,r10} obtained the multiple solutions problem
of $H$-surfaces, and Struwe \cite{r9} proved the existence of surfaces of constant mean curvature
$H$ with free boundaries.
For variable $H$, there are recent works by
 Caldiroli and Musina \cite{r13}and Duzzer and Grotowski \cite{r14}.
Caldiroli adn Musina  in
 \cite{r13}
considered   system  \eqref{e0}  with small boundary data, and proved
blowup phenomena and nonexistence results.
The existence of solutions to the
system for non-constant  $H$ in higher dimensional compact Riemannian manifolds
without boundary was proved by Duzzer-Grotowski \cite{r14}.



In this paper, we  study an  initial-boundary value  system for the
heat flow of the equation of $H$-surface:
\begin{equation}\label{e1}
\left\{
\begin{aligned}
u_{t}&=\Delta u-2 H(u) u_{x} \wedge u_{y},} {\text { in } \Omega \times(0, \infty), \\
\left.u\right|_{t=0}&=u_{0}, \text { in } \Omega, \\
\left.u\right|_{\partial \Omega}&=\chi,
\end{aligned}
\right.
\end{equation}
where
$u_{0} \in H^{1}(\Omega), \chi \in H^{\frac{1}{2}}(\partial \Omega), $  and $\left.u_{0}\right|_{\partial \Omega}=\chi$.
Throughout this paper, we assume that
\begin{equation}\label{e2}
H(u) \equiv  H\equiv const>0, \text { and } \chi=0.
\end{equation}

Struwe \cite{r9}, by the  assumption $|H(u)| \cdot|\chi|_{L^{\infty}(\partial \Omega)}<1$,
proved that the equations \eqref{e1} with  the condition
with free boundaries admit a unique
solution. By using the theorems and methods in \cite{r15, r9}, Rey \cite{r12} proved that
if $u_{0}(x) \in W^{1,2}(\Omega, \mathbb{R}^{3}) \cap L^{\infty}(\Omega, \mathbb{R}^{3})$
and $\left.u_{0}(x)\right|_{\partial \Omega}=\chi,$ then system \eqref{e1} has a unique global
regular solution $u \in C^{1+\frac{\alpha}{2}, 2+\alpha}(\Omega \times(0,+\infty), \mathbb{R}^{3})$
under the assumption $|H(u)| \cdot|\chi|_{L^{\infty}(\partial \Omega)}<1$.  Chen and Levine
\cite{r16} removed the assumption  $|H(u)| \cdot|\chi|_{L^{\infty}(\partial \Omega)}<1$,
 and obtained the existence of regular solution to system
\eqref{e1} but added the following assumption
\begin{equation}
\int_{\Omega}|\nabla u|^{2}(\cdot, t) \leq \int_{\Omega}|\nabla u|^{2}(\cdot, s), \quad \text { for } 0 \leq s \leq t,
\end{equation}
which is the main difference between the heat flow of the equation of $H$-surface and the
heat flow of harmonic maps.  The existence of weak solutions
and short-time regularity for the $H$-surface flow were considered by
 B\"{o}gelein, Duzar and Scheven \cite{r17,r18}. If $\chi \equiv 0,$ Huang, Tan and Wang \cite{r4} gave
sufficient conditions with  low initial energy   such that the heat flow develops finite time singularity.

In this article, we consider the heat flow system of
H-surface  with low initial energy, critical initial energy  and high initial energy.
The results in our paper will be  obtained by the modified   potential well method.
Potential well method, which was first put forward to consider semi-linear hyperbolic
initial boundary value problem by Payne and Sattinger \cite{r20,r21} around 1970s,
is a powerful tool in studying the long time behaviors of solutions of some
evolution equations. The potential well is   defined by the level set of energy
functional and  the derivative functional. It is generally true that solutions
starting inside the well are global in time, solutions starting outside the well
and at an unstable point blow up in finite time. After the pioneer work of
Sattinger and Payne, some authors \cite{r22,r23,r24,r5,r26,r19} used the method to study the global
  existence and nonexistence of solutions  for various nonlinear evolution equations with initial
   boundary value problem.  In \cite{r24,r5}, Liu  et al. modified and improved the method by introducing
    a family of potential wells which include the known potential well as a special case.
    The modified potential well method has been used to study semilinear pseudo-parabolic equations
    \cite{r19} and fourth-order parabolic equation \cite{r6}.  In this paper, we use the modified
    potential well method to obtain global existence and blow up in finite time of solutions when
    the initial energy is low, critical and high, respectively.  When the initial energy is low,
    similar results are obtained in \cite{r4}, but our result is more general, moreover, we prove
    a more precise decay rate of  $|u|_2$.

\section{Preliminaries}
Throughout this paper, we denote the $L^{2}(\Omega)$ norm,
$H_0^{1}(\Omega)$ norm by $|\cdot|_{2}$, $\|\cdot\|$, respectively.
And $(\cdot, \cdot)$ is used to denote the inner product
in $L^{2}(\Omega)$ . In order to state our main results precisely, let us introduce some notations and sets,
and then investigate their basic properties.

For $u\in H_0^1(\Omega)$, we set
\begin{equation}
\begin{aligned}
E(u) &=e(u)+V_{H}(u)=\frac{1}{2} \int_{\Omega}|\nabla u|^{2}+\frac{2}{3} \int_{\Omega} H(u) u \cdot u_{x} \wedge u_{y}, \\
 D(u) &=\int_{\Omega}|\nabla u|^{2}+2\int_{\Omega}H(u) uu_{x} \wedge u_{y}.
   \end{aligned}
\end{equation}
The Nehari manifold is defined by
\begin{equation}
\mathcal{N}=\left\{u \in H_{0}^{1}(\Omega) : D(u)=0, u \neq 0\right\},
\end{equation}
which can be  separated into the two unbounded sets
\begin{equation}
\begin{aligned}
\mathcal{N}_{+} &=\left\{u \in H_{0}^{1}(\Omega) : D(u)>0\right\}, \\
 \mathcal{N}_{-} &=\left\{u \in H_{0}^{1}(\Omega) : D(u)<0\right\}.
 \end{aligned}
\end{equation}
The  potential well and its corresponding set are defined respectively as
\begin{equation}
\begin{aligned}
  W &=\left\{u \in H_{0}^{1}(\Omega) : D(u)>0, E(u)<d\right\} \cup\{0\}, \\
   V &=\left\{u \in H_{0}^{1}(\Omega) : D(u)<0, E(u)<d\right\},
   \end{aligned}
\end{equation}
where,
$$ d=\min _{u \in H_{0}^{1}(\Omega) \backslash\{0\}} \max _{s \geqslant 0} E(s u)=\inf _{u \in \mathcal{N}} E(u),$$
is the depth of the potential well $W$.

Now let us define the level set
\begin{equation}
E^{\alpha}=\left\{u \in H_{0}^{1}(\Omega): E(u)<\alpha\right\}.
\end{equation}
Furthermore, by the definition of $E(u), \mathcal{N}, E^{\alpha}$ and $d,$
we easily know that
\begin{equation}
\mathcal{N}_{\alpha}=\mathcal{N} \cap E^{\alpha} \equiv\left\{u \in \mathcal{N} : \|u\|<\sqrt{6\alpha}\right\}
 \neq \varnothing \quad \text { for all } \alpha>d.
\end{equation}
We now define
\begin{equation}\label{ee123}
\lambda_{\alpha}=\inf \left\{|u|_2 : u \in \mathcal{N}_{\alpha}\right\}, \quad \Lambda_{\alpha}
=\sup \left\{|u|_2: u \in \mathcal{N}_{\alpha}\right\} \quad \text { for all } \alpha>d.
\end{equation}
It is clear that $\lambda_{\alpha}$ is nonincreasing and $\Lambda_{\alpha}$ is nondecreasing with respect to $\alpha$ .
We also introduce the following sets
\begin{equation}
\begin{aligned}
 \mathcal{B} &=\left\{u_{0} \in H_{0}^{1}(\Omega) : \text { the solution } u=u(t) \text { of }(1) \text { blows up in finite time }\right\}, \\
  \mathcal{G} &=\left\{u_{0} \in H_{0}^{1}(\Omega) : \text { the solution } u=u(t) \text { of }(1) \text { exists for all } t>0\right\}, \\
\mathcal{G}_{o} &=\left\{u_{0} \in G : u(t) \mapsto 0 \text { in } H_{0}^{1}(\Omega) \text { as } t \rightarrow \infty\right\}.
\end{aligned}
\end{equation}
For  $0<\delta<\frac{3}{2}$ , let us define the modified functional and Nehari manifold as follows:
\begin{equation}
\begin{aligned} D_{\delta}(u) &=\delta\|\nabla u\|_{2}^{2}+2\int_{\Omega}H(u)u\cdot u_x \wedge u_{y},\\
\mathcal{N}_{\delta} &=\left\{u \in H_{0}^{1}(\Omega) :  D_{\delta}(u)=0,\|u\| \neq 0\right\}, \\
d_{\delta} &=\inf _{u \in \mathcal{N}_{\delta}} E(u), \\
  r(\delta) &=\frac{2\sqrt{2\pi}\delta}{H}.
\end{aligned}
\end{equation}
Then we can  define the modified potential wells and their corresponding sets as follows:
\begin{equation}
\begin{aligned}
W_{\delta} &=\left\{u \in H_{0}^{1}(\Omega) : D_{\delta}(u)>0, E(u)<d(\delta)\right\} \cup\{0\}, \\
 V_{\delta} &=\left\{u \in H_{0}^{1}(\Omega) : D_{\delta}(u)<0, E(u)<d(\delta)\right\}, \\
  B_{\delta} &=\left\{u \in H_{0}^{1}(\Omega) :\|u\|<r(\delta)\right\}, \\
  B_{\delta}^{c} &=\left\{u \in H_{0}^{1}(\Omega) :\|u\|>r(\delta)\right\}.
  \end{aligned}
\end{equation}

For future convenience, we give some useful lemmas which will play
an important role in the proof of our main results. We  first recall
the following isoperimetric inequality, whose proof can be found in \cite{r1} and \cite{r2}.
\begin{lemma}[Isoperimetric inequality]\label{l11}
For any  $u \in H_{0}^{1}\left(\Omega ; \mathbb{R}^{3}\right)$, there holds
\begin{equation}
\int_{\Omega}|\nabla u|^{2} \geq \sqrt[3]{32 \pi}\left|\int_{\Omega} u \cdot u_{x} \wedge u_{y}\right|^{2 / 3}.
\end{equation}
\end{lemma}

\begin{lemma}\label{l25} Let $u_0\in H_0^1(\Omega)$.
\begin{itemize}\addtolength{\itemsep}{-1.5 em} \setlength{\itemsep}{-5pt}
\item[\emph{(1)}]
If $0<\|u\|<r(\delta),$ then $D_{\delta}(u)>0 .$ In particular, if $0<\|u\|_{H_{0}^{1}}<r(1),$ then $D(u)>0$;
\item[\emph{(2)}]  If $D_{\delta}(u)<0,$ then $\|u\|>r(\delta).$ In particular, if $D(u)<0,$ then $\|u\|>r(1)$;
\item[\emph{(3)}] If $D_{\delta}(u)=0,$ then $\|u\| \geqslant r(\delta)$ or $\|u\|=0 .$ In particular, if $D(u)=0,$ then $\|u\|
 \geqslant r(1)$
or $\|u\|=0$;
\item[\emph{(4)}] If $D_{\delta}(u)=0$ and $\|u\| \neq 0,$ then $E(u)>0$ for $0<\delta<\frac{3}{2}, E(u)=0$ for $\delta=\frac{3}{2}$,
$E(u)<0$ for $\delta>\frac{3}{2}$.
\end{itemize}
\end{lemma}

\begin{proof}
(1) Since $0<\|u\|<r(\delta)$, by the Isoperimetric inequality, we have
 \begin{equation}
 \left|\int_{\Omega} u \cdot u_{x} \wedge u_{y}\right|\leq \frac{1}{4\sqrt{2\pi}}\left(\int_{\Omega}|\nabla u|^{2}\right)^{\frac{3}{2}}.
\end{equation}
So from the assumption $0<\|u\|<r(\delta)=\frac{2\sqrt{2\pi}\delta}{H}$, we obtain
\begin{align}\label{e4}
   D_{\delta}(u) &=\delta\|\nabla u\|_{2}^{2}+2\int_{\Omega}H(u)u\cdot u_x \wedge u_{y}\nonumber\\
   &\geq\delta\|\nabla u\|_{2}^{2}-\frac{H}{2\sqrt{2\pi}}\left(\int_{\Omega}|\nabla u|^{2}\right)^{\frac{3}{2}}\nonumber\\
   &\geq\|\nabla u\|_{2}^{2}\left(\delta-\frac{H}{2\sqrt{2\pi}}\left(\int_{\Omega}|\nabla u|^{2}\right)^{\frac{1}{2}}\right)>0.
\end{align}

(2) By the assumption $D_{\delta}(u)<0$ and the Isoperimetric inequality, we have
\begin{align}\label{e5}
   0>D_{\delta}(u) &=\delta\|\nabla u\|_{2}^{2}+2\int_{\Omega}H(u)u\cdot u_x \wedge u_{y}\nonumber\\
   &\geq\delta\|\nabla u\|_{2}^{2}-\frac{H}{2\sqrt{2\pi}}\left(\int_{\Omega}|\nabla u|^{2}\right)^{\frac{3}{2}}\nonumber\\
   &\geq\|\nabla u\|_{2}^{2}\left(\delta-\frac{H}{2\sqrt{2\pi}}\left(\int_{\Omega}|\nabla u|^{2}\right)^{\frac{1}{2}}\right).
\end{align}
Hence, $\|u\|>r(\delta).$

(3) By the assumption $D_{\delta}(u)=0$ and the Isoperimetric inequality, we have
\begin{align}\label{e6}
   0=D_{\delta}(u) &=\delta\|\nabla u\|_{2}^{2}+2\int_{\Omega}H(u)u\cdot u_x \wedge u_{y}\nonumber\\
   &\geq\delta\|\nabla u\|_{2}^{2}-\frac{H}{2\sqrt{2\pi}}\left(\int_{\Omega}|\nabla u|^{2}\right)^{\frac{3}{2}}\nonumber\\
   &\geq\|\nabla u\|_{2}^{2}\left(\delta-\frac{H}{2\sqrt{2\pi}}\left(\int_{\Omega}|\nabla u|^{2}\right)^{\frac{1}{2}}\right).
\end{align}
Hence, $\|u\|\geq r(\delta)$ or $u=0$.

\noindent(4) We easily know that
\begin{align}\label{e7}
E(u)&=\frac{1}{2} \int_{\Omega}|\nabla u|^{2}+\frac{2}{3} \int_{\Omega} H(u) u \cdot u_{x} \wedge u_{y}\nonumber\\
&=\left(\frac{1}{2}-\frac{\delta}{3}\right)\|u\|^2+\frac{1}{3}D_{\delta}(u)\nonumber\\
&=\left(\frac{1}{2}-\frac{\delta}{3}\right)\|u\|^2.
\end{align}
Then  using \eqref{e7}, we can  prove the conclusion.
\end{proof}

\begin{lemma}\label{l26}
\begin{itemize}\addtolength{\itemsep}{-1.5 em} \setlength{\itemsep}{-5pt}
\item[\emph{(1)}] $d(\delta) \geqslant a(\delta) r^{2}(\delta)$ for $a(\delta)=\frac{1}{2}-\frac{\delta}{3}, 0<\delta<\frac{3}{2}$,
\item[\emph{(2)}] $\lim _{\delta \rightarrow 0} d(\delta)=0, d\left(\frac{3}{2}\right)=0$ and $d(\delta)<0$ for $\delta>\frac{3}{2}$,
\item[\emph{(3)}] $d(\delta)$ is increasing on $0<\delta \leqslant 1,$ decreasing on $1 \leqslant \delta \leqslant \frac{3}{2}$
and takes the maximum $d=$ $d(1)$ at $\delta=1.$
\end{itemize}
\end{lemma}
\begin{proof}
(1) If $u\in \mathcal{N}_{\delta}$,  by Lemma \ref{l25} (3), then $\|u\| \geqslant r(\delta)$. Moreover, we can deduce
\begin{align}\label{e8}
E(u)&=\frac{1}{2} \int_{\Omega}|\nabla u|^{2}+\frac{2}{3} \int_{\Omega} H(u) u \cdot u_{x} \wedge u_{y}\nonumber\\
&=\left(\frac{1}{2}-\frac{\delta}{3}\right)\|u\|^2+\frac{1}{3}D_{\delta}(u)\nonumber\\
&=\left(\frac{1}{2}-\frac{\delta}{3}\right)\|u\|^2\geq a(\delta) r^2(\delta).
\end{align}
Hence, $d(\delta) \geqslant a(\delta) r^{2}(\delta)$.

(2)
We easily know that
$$E(\lambda u)=\frac{\lambda^2}{2} \int_{\Omega}|\nabla u|^{2}+\frac{2\lambda^3}{3} \int_{\Omega} H(u) u \cdot u_{x} \wedge u_{y}.$$
Hence, \begin{equation}\label{e9}
\lim _{\lambda \rightarrow 0} E(\lambda u)=0.
\end{equation}
And if we let $\lambda u\in \mathcal{N}_{\delta}$, then $\lambda u$ satisfies
 $$0=D_{\delta}(\lambda u) =\delta \lambda^2\|\nabla u\|_{2}^{2}+2\lambda^3 \int_{\Omega}H(u)u\cdot u_x \wedge u_{y}.$$
Then, we obtain
\begin{equation}\label{e10}
 \lambda =\frac{\delta\|\nabla u\|_{2}^{2}}{2\int_{\Omega}H(u)u\cdot u_x \wedge u_{y}},
\end{equation}
which yields
\begin{equation}
\lim _{\delta \rightarrow 0} \lambda(\delta)=0.
\end{equation}
Now \eqref{e9} implies that
\begin{equation}
\lim _{\delta \rightarrow 0} E(\lambda u)=\lim _{\lambda \rightarrow 0} E(\lambda u)=0,
\end{equation}
 and
 \begin{equation}
\lim _{\delta \rightarrow 0} d(\delta)=0.
\end{equation}

It is easy to see that from \eqref{e8}
$$d\left(\frac{3}{2}\right)=0\ \  \mbox{and} \ \ d(\delta)<0 \ \ \mbox{for} \ \ \delta>\frac{3}{2}.$$

The proof is complete.

(3)  We need to prove that for any  $0<\delta^{\prime}<\delta^{\prime \prime}<1$ or $1<\delta^{\prime \prime}<\delta^{\prime}<\frac{3}{2}$
and for any $u \in \mathcal{N}_{\delta^{\prime \prime}},$ there is a $v \in \mathcal{N}_{\delta^{\prime}}$ and a constant $\varepsilon\left(\delta^{\prime}, \delta^{\prime \prime}\right)$ such that $E(v)<E(u)-\varepsilon\left(\delta^{\prime}, \delta^{\prime \prime}\right)$.
Indeed, by the definition of \eqref{e10}, we easily know that $D_{\delta}(\lambda(\delta) u)=0$  and $\lambda\left(\delta^{\prime \prime}\right)=1$.
Let $h(\lambda)=E(\lambda u)$, we have
\begin{equation}
\begin{aligned} \frac{d}{d \lambda} h(\lambda) &
=\frac{1}{\lambda}\left((1-\delta)\|\lambda u\|^{2}+D_{\delta}(\lambda u)\right) \\
&=(1-\delta) \lambda\|u\|^{2}.
 \end{aligned}
\end{equation}
Take $v=\lambda\left(\delta^{\prime}\right) u,$ then $v \in \mathcal{N}_{\delta^{\prime}}$.

 For $0<\delta^{\prime}<\delta^{\prime \prime}<1,$ we obtain
 \begin{equation}
\begin{aligned} E(u)-E(v) &=h(1)-h\left(\lambda\left(\delta^{\prime}\right)\right) \\
&>\left(1-\delta^{\prime \prime}\right) r^{2}\left(\delta^{\prime \prime}\right) \lambda\left(\delta^{\prime}\right)\left(1-\lambda\left(\delta^{\prime}\right)\right)
\equiv \varepsilon\left(\delta^{\prime}, \delta^{\prime \prime}\right).
\end{aligned}
\end{equation}
 For $1<\delta^{\prime \prime}<\delta^{\prime}<\frac{3}{2},$ we obtain
 \begin{equation}
\begin{aligned} E(u)-E(v) &=h(1)-h\left(\lambda\left(\delta^{\prime}\right)\right) \\
 &>\left(\delta^{\prime \prime}-1\right) r^{2}\left(\delta^{\prime \prime}\right) \lambda\left(\delta^{\prime \prime}\right)\left(\lambda\left(\delta^{\prime}\right)-1\right)
 \equiv \varepsilon\left(\delta^{\prime}, \delta^{\prime \prime}\right).
 \end{aligned}
\end{equation}
 Hence, the proof is complete.

\end{proof}

\begin{lemma}\label{l27} Let $u_{0} \in H_{0}^{1}(\Omega)$ and $0<\delta<\frac{3}{2} .$ If  $E(u) \leqslant d(\delta),$  then we have
\begin{itemize}\addtolength{\itemsep}{-1.5 em} \setlength{\itemsep}{-5pt}
\item[\emph{(1)}] If $D_{\delta}(u)>0,$ then $\|u\|^{2}<\frac{d(\delta)}{a(\delta)},$
 where $a(\delta)=\frac{1}{2}-\frac{\delta}{3} .$ In particular, if $D(u) \leqslant d$ and $D(u)>0,$ then
    \begin{equation}
\|u\|^{2}<6d.
\end{equation}
\item[\emph{(2)}] If $\|u\|^{2}>\frac{d(\delta)}{a(\delta)},$ then $D_{\delta}(u)<0 .$ In particular, if $E(u) \leqslant d$ and
\begin{equation}
\|u\|^{2}>6d,
\end{equation}
then $D(u)<0$.
\item[\emph{(3)}]  If $D_{\delta}(u)=0,$ then $\|u\|^{2} \leqslant \frac{d(\delta)}{a(\delta)} .$ In particular, if $E(u) \leqslant d$ and $D(u)=0,$ then  \begin{equation}
\|u\|^{2} \leqslant 6d.
\end{equation}
\end{itemize}
\end{lemma}
\begin{proof}
(1)For  $0<\delta<\frac{3}{2}$, we see that
\begin{align}\label{e11}
E(u)&=\frac{1}{2} \int_{\Omega}|\nabla u|^{2}+\frac{2}{3} \int_{\Omega} H(u) u \cdot u_{x} \wedge u_{y}\nonumber\\
&=\left(\frac{1}{2}-\frac{\delta}{3}\right)\|u\|^2+\frac{1}{3}D_{\delta}(u)\nonumber\\
&=a(\delta)\|u\|^2\leq d(\delta).
\end{align}
Therefore,  $$\|u\|^{2}<\frac{d(\delta)}{a(\delta)}.$$

Finally, (2) and (3) follow from (\ref{e11}).
\end{proof}

\begin{lemma}\label{l2801}
Let $u \in H_{0}^{1}(\Omega)$. We have
\begin{itemize}\addtolength{\itemsep}{-1.5 em} \setlength{\itemsep}{-5pt}
\item[\emph{(1)}] 0 is away from both $N$ and $\mathcal{N}_{-},$ i.e. dist $(0, \mathcal{N})>0,$  $dist \left(0, N_{-}\right)>0$,
\item[\emph{(2)}] For any $\alpha>0,$ the set $E^{\alpha} \cap \mathcal{N}_{+}$ is bounded in $H_{0}^{1}(\Omega)$.
\end{itemize}
\end{lemma}
\begin{proof}
(1) If $u\in \mathcal{N}$, then we have
\begin{equation}
\begin{aligned}
d & \leq  E(u)=\frac{1}{2} \int_{\Omega}|\nabla u|^{2}+\frac{2}{3} \int_{\Omega} H(u) u \cdot u_{x} \wedge u_{y} \\
&= \frac{1}{6}\int_{\Omega}|\nabla u|^{2}+\frac{1}{3}D(u)= \frac{1}{6}\int_{\Omega}|\nabla u|^{2}.
 \end{aligned}
\end{equation}
If  $u\in \mathcal{N}_{-}$, then we have
\begin{equation}
\begin{aligned}
d & \leq  E(u)=\frac{1}{2} \int_{\Omega}|\nabla u|^{2}+\frac{2}{3} \int_{\Omega} H(u) u \cdot u_{x} \wedge u_{y} \\
&= \frac{1}{6}\int_{\Omega}|\nabla u|^{2}+\frac{1}{3}D(u)\leq  \frac{1}{6}\int_{\Omega}|\nabla u|^{2}.
 \end{aligned}
\end{equation}
Hence, 0 is away from both $N$ and $\mathcal{N}_{-},$ i.e. dist $(0, \mathcal{N})>0,$  $dist \left(0, N_{-}\right)>0$.

(2) Since $E(u)<\alpha$ and $D(u)>0,$ we obtain
\begin{equation}
\begin{aligned}
\alpha & >  E(u)=\frac{1}{2} \int_{\Omega}|\nabla u|^{2}+\frac{2}{3} \int_{\Omega} H(u) u \cdot u_{x} \wedge u_{y} \\
&= \frac{1}{6}\int_{\Omega}|\nabla u|^{2}+\frac{1}{3}D(u)> \frac{1}{6}\int_{\Omega}|\nabla u|^{2}.
 \end{aligned}
\end{equation}
Hence,  for any $\alpha>0,$ the set $E^{\alpha} \cap \mathcal{N}_{+}$ is bounded in $H_{0}^{1}(\Omega)$.
\end{proof}

\section{Low initial energy $E(u_0)<d$.}
The goal of this section is to prove Theorem \ref{t21}--\ref{t23}.
 A threshold result for the global solutions and finite time blowup  will be given.


\begin{theorem}\label{t20}
Assume that $u_{0} \in H_{0}^{1}(\Omega)$,
$T$ is the maximal existence time of $u$, and
$0<e<d, \delta_{1}<\delta_{2}$
are  two roots of equation $d(\delta)=e$.  We have
\begin{itemize}\addtolength{\itemsep}{-1.5 em} \setlength{\itemsep}{-5pt}
\item[\emph{(1)}]If $D(u_0) > 0$, all weak solutions $u$ of system  (\ref{e1})
with $E\left(u_{0}\right)=e$ belong to $W_{\delta}$ for $\delta_{1}<\delta<\delta_{2}, 0 \leqslant t<T$.
\item[\emph{(2)}] If $D(u_0) < 0$, all weak solutions $u$ of system  (\ref{e1})
with $E\left(u_{0}\right)=e$ belong to $V_{\delta}$ for $\delta_{1}<\delta<\delta_{2}, 0 \leqslant t<T$.

\end{itemize}
\end{theorem}

\begin{theorem}\label{t21} (Global existence)
 Assume that
$u_{0} \in H_{0}^{1}(\Omega)$, $E\left(u_{0}\right)<d$, $D\left(u_{0}\right)>0$. Then system (\ref{e1})
has a global  solution
$u(t) \in L^{\infty}(0, \infty ; H_{0}^{1}(\Omega))$
 and $u(t) \in W$ for $0 \leqslant t<\infty .$
\end{theorem}

\begin{remark}\label{rema1}
Result similar to Theorem \ref{t21} is obtained in  \cite{r4}. But our proof is different to \cite{r4}. In fact,
using  the modified potential well method we can obtain the more general conclusion:

 If  the assumption $D\left(u_{0}\right)>0$ is replaced by $D_{\delta_{2}}\left(u_{0}\right)>0,$ where $\delta_{1}<\delta_{2}$
are the two roots of equation $d(\delta)=E\left(u_{0}\right),$  then system  \eqref{e1} admits a global weak solution.
\end{remark}

The following result is obtained in  \cite{r4}.
But our proof is different from the proof in \cite{r4}.
 For the reader's
convenience, we will give the detailed proof.
\begin{theorem}\label{t22}
Assume that   $u_{0} \in H_{0}^{1}(\Omega),$  $E\left(u_{0}\right)<d$
and $D\left(u_{0}\right)<0.$ Then the weak solution $u(t)$ of system  \eqref{e1}
blows up in finite time, that is, there exists $a$
$T>0$ such that
$$
\lim _{t \rightarrow T} \int_{0}^{t}|u(\tau)|_2 d \tau=+\infty
$$
\end{theorem}
\begin{remark}\label{rema3}
Assume that
$u_{0} \in H_{0}^{1}(\Omega)$,$E\left(u_{0}\right)<d .$ When $D\left(u_{0}\right)>0,$  system  (\ref{e1}) has a global solution.
When $D\left(u_{0}\right)<0,$ system  \eqref{e1} does not admit any global weak solution.
\end{remark}

\begin{theorem}\label{t23}
Assume that
$u_{0} \in H_{0}^{1}(\Omega)$, $E\left(u_{0}\right)<d$ and $D\left(u_{0}\right)>0,$  $\delta_{1}<\delta_{2}$ are the two roots of equation
$d(\delta)=E(u_0).$  Then, for the global weak solution $u$ of system  \eqref{e1}, it holds
\begin{equation}
|u|_2^2 \leqslant\left|u_{0}\right|_2^{2} e^{-2\left(1-\delta_{1}\right) t}, \quad 0 \leqslant t<\infty.
\end{equation}
\end{theorem}
\begin{remark}\label{rema4}
 In comparison with  the  decay rate of  \cite{r4},  the   result of the decay rate of  $|u|_2$ in Theorem \ref{t23} is much more  precise.
\end{remark}

In order to prove Theorem \ref{t31} -- \ref{t33}, we need the following lemmas:

\begin{lemma}\label{l29}
  For $0<T \leq \infty,$ assume that $u : \Omega \times[0, T) \rightarrow \mathbb{R}^{3}$ is a weak
  solution to system  \eqref{e1}. Then it holds
  \begin{equation}\label{e12}
\int_{t_{1}}^{t_{2}} \int_{\Omega}\left|u_{t}\right|^{2}+E\left(u\left(t_{2}\right)\right)
=E\left(u\left(t_{1}\right)\right), \quad \forall  t_{1},t_{2}\in (0,T).
\end{equation}
  \begin{proof}
   Multiplying \eqref{e1} by $u_t$ and integrating over $\Omega$ via the integration by
parts we get \eqref{e12}.
  \end{proof}

\end{lemma}

\begin{lemma}\label{l28}
If $0<E(u)<d$ for some $u \in H_{0}^{1}(\Omega),$ and $\delta_{1}<1<\delta_{2}$ are the two roots of equation
$d(\delta)=E(u),$ then  the sign of $D_{\delta}(u)$ doesn't change for $\delta_{1}<\delta<\delta_{2}$ .
\end{lemma}

\begin{proof}
Since $E(u)>0$, we have  $\|u\| \neq 0.$
If the sign of $D_{\delta}(u)$ is changeable for $\delta_{1}<\delta<\delta_{2},$ then
we choose $\bar{\delta} \in\left(\delta_{1}, \delta_{2}\right)$ such that $D_{\bar{\delta}}(u)=0 .$
Hence, by the definition of $d(\bar{\delta})$, we can obtain $E(u) \geqslant d(\bar{\delta}),$ which contradicts
$E(u)=d\left(\delta_{1}\right)=d\left(\delta_{2}\right)<d(\bar{\delta})$ (by Lemma \ref{l26} (3)).
\end{proof}

\begin{definition}\label{d21}
(Maximal existence time). Assume that  $u(t)$ is  a weak solution of system  (\ref{e1}).
The maximal existence time $T$ of $u(t)$  is defined as follows:
\begin{itemize}\addtolength{\itemsep}{-1.5 em} \setlength{\itemsep}{-5pt}
\item[\emph{(1)}] If $u(t)$ exists for $0 \leqslant t<\infty,$ then $T=+\infty$.
\item[\emph{(2)}] If there is a $t_{0} \in(0, \infty)$ such that $u(t)$ exists for $0 \leqslant t<t_{0},$ but doesn't exist at $t=t_{0},$
then $T=t_{0}.$
\end{itemize}
\end{definition}

 \begin{proof}[\textbf{Proof of theorem \ref{t20}}]
(1) Let $u(t)$ be any weak solution of system  (\ref{e1}) with $E\left(u_{0}\right)=e,$ $D\left(u_{0}\right)>0,$ and $T$ be the
maximal existence time of $u(t).$ Using $E\left(u_{0}\right)=e, D\left(u_{0}\right)>0$ and Lemma \ref{l28},
 we have  $D_{\delta}\left(u_{0}\right)>0$
  and $E\left(u_{0}\right)<d(\delta) .$ So $u_{0}(x) \in W_{\delta}$ for $\delta_{1}<\delta<\delta_{2}.$
  We need to prove that  $u(t) \in W_{\delta}$  for  $\delta_{1}<\delta<\delta_{2}$  and
$0<t<T.$   Indeed, if this is not the conclusion, from time continuity of $D(u)$ we assume that there must
exist a $\delta_{0} \in\left(\delta_{1}, \delta_{2}\right)$ and $t_{0} \in(0, T)$ such that
$u\left(t_{0}\right) \in \partial W_{\delta_{0}},$ and
$D_{\delta_{0}}\left(u\left(t_{0}\right)\right)=0,\left\|u\left(t_{0}\right)\right\| \neq 0$ or
$E\left(u\left(t_{0}\right)\right)=d\left(\delta_{0}\right).$ From the energy  equality
   \begin{equation}\label{e13}
\int_{0}^{t} \int_{\Omega}\left|u_{\tau}\right|^{2}+E\left(u\left(t\right)\right)
=E(u_0)<d(\delta), \ \delta_{1}<\delta<\delta_{2}, \quad 0 \leqslant t<T,
\end{equation}
  we easily know that $E\left(u\left(t_{0}\right)\right) \neq d\left(\delta_{0}\right).$ If $D_{\delta_{0}}\left(u\left(t_{0}\right)\right)=0,\left\|u\left(t_{0}\right)\right\| \neq 0,$
  then by the definition of $d(\delta)$
we obtain  $E\left(u\left(t_{0}\right)\right) \geqslant d\left(\delta_{0}\right),$ which contradicts $(\ref{e13}).$

  (2) Let $u(t)$ be any weak solution of system  (\ref{e1}) with $E\left(u_{0}\right)=e,$ $D\left(u_{0}\right)<0,$ and $T$ be the
maximal existence time of $u(t).$ Using $E\left(u_{0}\right)=e, D\left(u_{0}\right)<0$ and Lemma \ref{l28},
 we have  $D_{\delta}\left(u_{0}\right)<0$
  and $E\left(u_{0}\right)<d(\delta) .$ So $u_{0} \in V_{\delta}$ for $\delta_{1}<\delta<\delta_{2}.$
  We need to prove that  $u(t) \in V_{\delta}$  for  $\delta_{1}<\delta<\delta_{2}$  and
$0<t<T.$   Indeed, if this is not the conclusion, from time continuity of $D(u)$ we assume that there must
exist a $\delta_{0} \in\left(\delta_{1}, \delta_{2}\right)$ and $t_{0} \in(0, T)$ such that
$u\left(t_{0}\right) \in \partial V_{\delta_{0}},$ and
$D_{\delta_{0}}\left(u\left(t_{0}\right)\right)=0,$ or
$E\left(u\left(t_{0}\right)\right)=d\left(\delta_{0}\right).$  From the energy  equality \eqref{e13},
  we easily know that $E\left(u\left(t_{0}\right)\right) \neq d\left(\delta_{0}\right).$
  If $D_{\delta_{0}}\left(u\left(t_{0}\right)\right)=0,$
  and $t_{0}$ is the first time such that
$D_{\delta_{0}}(u(t))=0,$ then $D_{\delta_{0}}(u(t))<0$ for $0 \leqslant t<T$.
 By Lemma \eqref{l25} (2), we have $\left\|u\left(t_{0}\right)\right\|>r\left(\delta_{0}\right)$ for $0 \leqslant t<T$.
 So, $\left\|u\left(t_{0}\right)\right\|>r\left(\delta_{0}\right)$ and
$E\left(u\left(t_{0}\right)\right) \neq d\left(\delta_{0}\right),$ which contradicts \eqref{e13}.

  \end{proof}

 \begin{proof}[\textbf{Proof of theorem \ref{t21}}]
From the standard argument  in \cite{r3}, we can prove the local existence result of \eqref{e1}
in a more general case of initial value $u_{0} \in H_0^{1}(\Omega)$ and
$u \in C^{0}\left(\left[0, T_{0}\right], H_0^{1}(\Omega)\right)$.

Using $E\left(u_{0}\right)<d, D\left(u_{0}\right)>0$ and Lemma \ref{l28},
 we have  $D_{\delta}\left(u_{0}\right)>0$
  and $E\left(u_{0}\right)<d(\delta).$ So $u_{0}(x) \in W_{\delta}$ for $\delta_{1}<\delta<\delta_{2}.$
  We need to prove that  $u(t) \in W_{\delta}$  for  $\delta_{1}<\delta<\delta_{2}$  and
$0<t<T.$   Indeed, if this is not the conclusion, from time continuity of $D(u)$ we assume that there must
exist a $\delta_{0} \in\left(\delta_{1}, \delta_{2}\right)$ and $t_{0} \in(0, T)$ such that
$u\left(t_{0}\right) \in \partial W_{\delta_{0}},$ and
$D_{\delta_{0}}\left(u\left(t_{0}\right)\right)=0,\left\|u\left(t_{0}\right)\right\| \neq 0$ or
$E\left(u\left(t_{0}\right)\right)=d\left(\delta_{0}\right).$ From the energy  equality
   \begin{equation}\label{e15}
\int_{0}^{t} \int_{\Omega}\left|u_{\tau}\right|^{2}+E\left(u\left(t\right)\right)
=E(u_0)<d(\delta), \ \delta_{1}<\delta<\delta_{2}, \quad 0 \leqslant t<T,
\end{equation}
  we easily know that $E\left(u\left(t_{0}\right)\right) \neq d\left(\delta_{0}\right).$ If $D_{\delta_{0}}\left(u\left(t_{0}\right)\right)=0,\left\|u\left(t_{0}\right)\right\| \neq 0,$
  then by the definition of $d(\delta)$
we obtain  $E\left(u\left(t_{0}\right)\right) \geqslant d\left(\delta_{0}\right),$ which contradicts $(\ref{e13}).$

\end{proof}

\begin{remark}\label{rema2}
If in Theorem \ref{t21} the condition $D_{\delta_{2}}\left(u_{0}\right)>0$ is replaced by $\left\|u_{0}\right\|<r\left(\delta_{2}\right),$
then system  \eqref{e1} has a global weak solution
$u(t) \in L^{\infty}(0, \infty ; H_{0}^{1}(\Omega))$ with $u_{t}(t) \in L^{2}(0, \infty ; H_{0}^{1}(\Omega))$ and the following result holds
\begin{equation}
\|u\|<\frac{d(\delta)}{a(\delta)}, \quad \delta_{1}<\delta<\delta_{2}, 0 \leqslant t<\infty,
\end{equation}
\begin{equation}
\int_{0}^{t}\left|u_{\tau}\right|^{2} d \tau<d(\delta), \quad \delta_{1}<\delta<\delta_{2}, 0 \leqslant t<\infty.
\end{equation}
In particular
\begin{equation}
\|u\|^{2}<\frac{d\left(\delta_{1}\right)}{a\left(\delta_{1}\right)},
\end{equation}
\begin{equation}
\int_{0}^{t}\left|u_{\tau}\right|^{2} d \tau<d\left(\delta_{1}\right), \quad 0 \leqslant t<\infty.
\end{equation}
\end{remark}

 \begin{proof}[\textbf{Proof of theorem \ref{t22}}] We argue by contradiction. Suppose that there would exist
a global weak solution $u(t)$. Set
\begin{equation}
f(t)=\int_{0}^{t} \int_{\Omega}|u|^{2}, t>0.
\end{equation}
Multiplying \eqref{e1} by $u$ and integrating over $\Omega \times(0, t),$ we get
\begin{equation}
\int_{\Omega}|u(t)|^{2}-\int_{\Omega}\left|u_{0}\right|^{2}
=-2 \int_{0}^{t} \int_{\Omega}\left(|\nabla u|^{2}+2 H(u) u \cdot u_{x} \wedge u_{y}\right).
\end{equation}
According to the definition of $f(t),$ we have $f^{\prime}(t)=\int_{\Omega}|u(t)|^{2}$ and hence
\begin{equation}
f^{\prime}(t)= \int_{\Omega}|u|^{2}=\int_{\Omega}\left|u_{0}\right|^{2}-2 \int_{0}^{t} \int_{\Omega}\left(|\nabla u|^{2}
+2 H(u) u \cdot u_{x} \wedge u_{y}\right),
\end{equation}
and
\begin{equation}\label{e17}
f^{\prime \prime}(t)=-2 \int_{\Omega}\left(|\nabla u|^{2}+2 H(u) u \cdot u_{x} \wedge u_{y}\right)=-2D(u).
\end{equation}
Now using  \eqref{e12},\eqref{e17}  and
\begin{align}\label{e16}
E(u)&=\frac{1}{2} \int_{\Omega}|\nabla u|^{2}+\frac{2}{3} \int_{\Omega} H(u) u \cdot u_{x} \wedge u_{y}\nonumber\\
&=\frac{1}{6}\|u\|^2+\frac{1}{3}D(u),
\end{align}
we can obtain

\begin{align}
f^{\prime \prime}(t) &= 6 \int_{0}^{t} \int_{\Omega}u^2_{\tau}d\tau
 +f^{\prime}(t)-6 E\left(u_{0}\right)\nonumber\\
 &=6 \int_{0}^{t} \int_{\Omega}u^2_{\tau}d\tau
 + \int_{\Omega}|u|^{2}-6 E\left(u_{0}\right).
\end{align}
Note that

\begin{align}
f(t)f^{\prime \prime}(t) &= f(t)\left[6 \int_{0}^{t} \int_{\Omega}u^2_{\tau}d\tau
 +f^{\prime}(t)-6 E\left(u_{0}\right)\right]\nonumber\\
 &=6\int_{0}^{t} \int_{\Omega}|u|^{2} \int_{0}^{t} \int_{\Omega}u^2_{\tau}d\tau
 +f(t)f^{\prime}(t)-6 E(u_{0})\int_{0}^{t} \int_{\Omega}|u|^{2}.
\end{align}
Hence, we have

\begin{align}
f(t)f^{\prime \prime}(t) -\frac{3}{2}(f^{\prime}(t))^2&
 \geq 6\int_{0}^{t} \int_{\Omega}|u|^{2} \int_{0}^{t} \int_{\Omega}u^2_{\tau}d\tau-6\left(\int_{0}^{t}\int_{\Omega}u_{\tau}\cdot u \mathrm{d} \tau\right)^{2}\nonumber\\
 &\ \ +f(t)f^{\prime}(t) -3f^{\prime}(t)\int_{\Omega}u_0^2- 6 E(u_{0})\int_{0}^{t} \int_{\Omega}|u|^{2}.
\end{align}
Making use of the Schwartz inequality, we have
\begin{align}\label{ee1}
f(t)f^{\prime \prime}(t) -\frac{3}{2}(f^{\prime}(t))^2&
 \geq f(t)f^{\prime}(t) -3f^{\prime}(t)\int_{\Omega}u_0^2- 6 E(u_{0})\int_{0}^{t} \int_{\Omega}|u|^{2}.
\end{align}

Next, we  distinguish two case:

(1) If $E\left(u_{0}\right) \leqslant 0,$ then
\begin{align}
f(t)f^{\prime \prime}(t) -\frac{3}{2}(f^{\prime}(t))^2&
 \geq f(t)f^{\prime}(t) -3f^{\prime}(t)\int_{\Omega}u_0^2.
\end{align}
Now we prove $D(u)<0$ for $t>0 .$  If not, we must be allowed to choose a $t_{0}>0$ such
that $D\left(u\left(t_{0}\right)\right)=0$ and $D(u)<0$ for $0 \leqslant t<t_{0} .$ From Lemma \ref{l25} (2),
we have $\|u\|>r(1)$  for
$0 \leqslant t<t_{0},\left\|u\left(t_{0}\right)\right\| \geqslant r(1)$ and $E\left(u\left(t_{0}\right)\right) \geqslant d,$
which contradicts \eqref{e13}.  From \eqref{e17} we have
$f^{\prime}(t)>0$ for $t \geqslant 0 .$ From $f^{\prime}(0)=\int_{\Omega}|u_0|^{2} \geqslant 0,$ we can know that there exists a $t_{0} \geqslant 0$ such
that $f^{\prime}\left(t_{0}\right)>0 .$ For $t \geqslant t_{0}$ we have
\begin{equation}
f(t) \geqslant f^{\prime}\left(t_{0}\right)\left(t-t_{0}\right)>f^{\prime}(0)\left(t-t_{0}\right).
\end{equation}
Hence, for sufficiently large $t$ , we  obtain
\begin{equation}
 f(t)>3\int_{\Omega}|u_0|^{2},
\end{equation}
then
$$f(t)f^{\prime \prime}(t) -\frac{3}{2}(f^{\prime}(t))^2>0.$$

(2) If $0<E\left(u_{0}\right)<d,$ then by Theorem \ref{t20} we have
$u(t) \in V_{\delta}$ for $1<\delta<\delta_{2}, t \geqslant 0,$ and $D_{\delta}(u)<0$,
$\|u\|>r(\delta)$ for $1<\delta<\delta_{2}, t \geqslant 0,$ where $\delta_{2}$ is the larger root of equation $d(\delta)=E\left(u_{0}\right) .$
Hence, $D_{\delta_{2}}(u) \leqslant 0$ and $\|u\| \geqslant r\left(\delta_{2}\right)$ for $t \geqslant 0.$
By (\ref{e17}), we have

\begin{equation}
\begin{aligned}
 f^{\prime \prime}(t) &=-2 D(u)=2\left(\delta_{2}-1\right)\|\nabla u\|_{2}^{2}-2 D_{\delta_{2}}(u), \\
 & \geqslant 2\left(\delta_{2}-1\right)|\nabla u|_{2}=2\left(\delta_{2}-1\right)\|u\|^{2} \geqslant 2\left(\delta_{2}-1\right) r^{2}\left(\delta_{2}\right), \quad t \geqslant 0, \\
  f^{\prime}(t) & \geqslant 2\left(\delta_{2}-1\right) r^{2}\left(\delta_{2}\right)t+f^{\prime}(0) \geqslant 2\left(\delta_{2}-1\right) r^{2}\left(\delta_{2}\right) t, \quad t \geqslant 0,
 \\ f(t) &\geqslant \left(\delta_{2}-1\right) r^{2}\left(\delta_{2}\right) t^{2}, \quad t \geqslant 0.
 \end{aligned}
\end{equation}
Therefore, for sufficiently large $t$, we infer
\begin{equation}
\frac{1}{2} f(t)>3\int_{\Omega}|u_0|^{2},\ \  \frac{1}{2} f^{\prime}(t) >6 E\left(u_{0}\right).
\end{equation}
Then,  (\ref{ee1}) implies that
\begin{equation}\label{ee2}
\begin{aligned}
f(t)f^{\prime \prime}(t) -\frac{3}{2}(f^{\prime}(t))^2&
 \geq f(t)f^{\prime}(t) -3f^{\prime}(t)\int_{\Omega}u_0^2- 6 E(u_{0})f(t).\nonumber\\
=&\left(\frac{1}{2} f(t)-3\int_{\Omega}|u_0|^{2}\right) f^{\prime}(t) \nonumber\\
 &\ \ \ +\left(\frac{1}{2} f^{\prime}(t)-6 E\left(u_{0}\right)\right) f(t)>0.
\end{aligned}
\end{equation}
The remainder of the proof is the same as that in \cite{r5}.
 \end{proof}

\begin{proof}[\textbf{Proof of theorem \ref{t23}}]
Multiplying \eqref{e1} by $v$, $v \in L^{\infty}\left(0, \infty; H_{0}^{1}(\Omega)\right)$, we have

\begin{equation}\label{e212}
 (u_t, v)+(\nabla u, \nabla v)=2H(u_x\wedge u_y, v).
\end{equation}
Letting $v = u$, (\ref{e212}) implies that
\begin{equation}\label{e213}
\frac{1}{2} \frac{d}{d t}\left|u\right|_2^{2}+D(u)=0, \quad 0 \leqslant t<\infty.
\end{equation}

From $0<E\left(u_{0}\right)<d, D\left(u_{0}\right)>0$ and Lemma \ref{t20},
we have $u(t) \in W_{\delta}$ for $\delta_{1}<\delta<\delta_{2}$ and
 $0 \leqslant$ $t<\infty,$ where $\delta_{1}<\delta_{2}$ are the two roots of equation
  $d(\delta)=E\left(u_{0}\right).$ Hence, we obtain  $D_{\delta}(u) \geqslant 0$ for
$\delta_{1}<\delta<\delta_{2}$ and $D_{\delta_{1}}(u) \geqslant 0$ for $0 \leqslant t<\infty .$ So, \eqref{e213} gives

\begin{equation}
\frac{1}{2} \frac{d}{d t}\left|u\right|_2^{2}+(1-\delta_1)\left|u\right|_2^{2}+D_{\delta}(u)=0, \quad 0 \leqslant t<\infty.
\end{equation}
Now  \eqref{e213}  implies that
\begin{equation}
\frac{1}{2} \frac{d}{d t}\left|u\right|_2^{2}+(1-\delta_1)\left|u\right|_2^{2}\leq 0, \quad 0 \leqslant t<\infty.
\end{equation}
and
\begin{equation}
|u|_2^{2} \leqslant\left|u_{0}\right|_2^{2}-2\left(1-\delta_{1}\right) \int_{0}^{t}|u(\tau)|^{2} d \tau, \quad 0 \leqslant t<\infty.
\end{equation}
By Gronwall's inequality, we have

\begin{equation}
|u|_{2}^{2} \leqslant\left|u_{0}\right|_{2}^{2} e^{-2\left(1-\delta_{1}\right) t}, \quad 0 \leqslant t<\infty.
\end{equation}
\end{proof}

\section{Critical initial energy $E(u_0)=d$.}
The goal of this section is to prove Theorem \ref{t31}--\ref{t33}.
\begin{theorem}\label{t31}
(Global existence) Assume
that $u_{0} \in H_{0}^{1}(\Omega)$,$E\left(u_{0}\right)=d$ and $D\left(u_{0}\right) \geqslant 0.$
Then system (\ref{e1}) has $a$ global weak solution $u(t) \in$
$L^{\infty}(0, \infty ; H_{0}^{1}(\Omega))$
and $u(t) \in \overline{W}=W \cup \partial W$ for $0 \leqslant t<\infty$
\end{theorem}

\begin{lemma}\label{l31}Assume that  $u\in H_0^1(\Omega)$, $\|\nabla u\|_{2}^{2} \neq 0$,  and $D(u)\geq 0$.  Then:
\begin{itemize}\addtolength{\itemsep}{-1.5 em} \setlength{\itemsep}{-5pt}
\item[\emph{(1)}]  $\lim _{\lambda \rightarrow 0} E(\lambda u)=0, \lim _{\lambda \rightarrow+\infty} E(\lambda u)=-\infty$,
\item[\emph{(2)}] On the interval $0<\lambda<\infty,$ there exists a unique $\lambda^{*}=\lambda^{*}(u),$ such that
\begin{equation}
\frac{d}{d \lambda} E\left.(\lambda u)\right|_{\lambda=\lambda^{*}}=0,
\end{equation}
\item[\emph{(3)}] $E(\lambda u)$ is increasing on $0 \leqslant \lambda \leqslant \lambda^{*},$ decreasing on $\lambda^{*} \leqslant \lambda<\infty$ and takes the maximum at $\lambda=\lambda^{*},$
    \item[\emph{(4)}] $D(\lambda u)>0$ for $0<\lambda<\lambda^{*}, D(\lambda u)<0$ for $\lambda^{*}<\lambda<\infty,$ and $D\left(\lambda^{*} u\right)=0$.
\end{itemize}
\end{lemma}
\begin{proof}
(1) Firstly, the assumption $D(u)\geq 0$ implies that  $$\int_{\Omega}H(u) u\cdot u_{x} \wedge u_{y}<0.$$
From the definition of $E(u)$, i.e.
$$E(u) =\frac{1}{2} \int_{\Omega}|\nabla u|^{2}+\frac{2}{3} \int_{\Omega} H(u) u \cdot u_{x} \wedge u_{y}, $$
and we see that

$$E(\lambda u) =\frac{\lambda^2}{2} \int_{\Omega}|\nabla u|^{2}+\frac{2\lambda^3}{3} \int_{\Omega} H(u) u \cdot u_{x} \wedge u_{y}.$$
Hence, we have
\begin{equation}
\lim _{ \lambda \rightarrow 0} E(\lambda u)=0\ \ \ \mbox{and}\ \ \ \lim _{\lambda \rightarrow+\infty} E(\lambda u)=-\infty.
\end{equation}

(2)  It is  easy to show that
$$\frac{d}{d\lambda}E(\lambda u) =\lambda\int_{\Omega}|\nabla u|^{2}+2\lambda^2\int_{\Omega} H(u) u \cdot u_{x} \wedge u_{y},$$
which leads to the conclusion.

(3) By Lemma \ref{l31} (2),  one has

\begin{equation}
\begin{aligned}
{\frac{d}{d \lambda} E(\lambda u)>0} & {\text { for } 0<\lambda<\lambda^{*}}, \\
 {\frac{d}{d \lambda} E(\lambda u)<0} & {\text { for } \lambda^{*}<\lambda<\infty},
\end{aligned}
\end{equation}
which leads to the conclusion.

(4)  The conclusion follows from
$$D(\lambda)=\frac{d}{d\lambda}E(\lambda u) =\lambda\int_{\Omega}|\nabla u|^{2}+2\lambda^2\int_{\Omega} H(u) u \cdot u_{x} \wedge u_{y}.$$

\end{proof}

\begin{proof}[\textbf{Proof of theorem \ref{t31}}]
 Firstly,  $E\left(u_{0}\right)=d$ implies that $\left\|u_{0}\right\|_{H_{0}^{1}} \neq 0.$
 Choose a sequence $\left\{\lambda_{m}\right\}$ such that $0<\lambda_{m}<1,$
$m=1,2, \ldots$ and $\lambda_{m} \rightarrow 1$ as $m \rightarrow \infty .$
 Let $u_{0 m}=\lambda_{m} u_{0}.$  We consider the following initial problem
\begin{equation}\label{e31}
\left\{
\begin{aligned}
u_{t}&=\Delta u-2 H(u) u_{x} \wedge u_{y},} {\text { in } \Omega \times(0, \infty), \\
\left.u\right|_{t=0}&=u_{0m}, \text { in } \Omega, \\
\left.u\right|_{\partial \Omega}&=\chi,  {t>0},
\end{aligned}
\right.
\end{equation}
 From $D\left(u_{0}\right) \geqslant 0$ and Lemma \ref{l31},  we have $\lambda^{*}=$
$\lambda^{*}\left(u_{0}\right) \geqslant 1.$ Thus, we get $D\left(u_{0 m}\right)=D\left(\lambda_{m} u_{0}\right)>0$
and $E\left(u_{0 m}\right)=E\left(\lambda_{m} u_{0}\right)<E\left(u_{0}\right)=d .$ From
Theorem \ref{t21}, it follows that for each $m$ problem \eqref{e31} admits a global weak solution $u_{m}(t) \in$
$L^{\infty}(0, \infty ; H_{0}^{1}(\Omega))$
with $u_{m t}(t) \in L^{2}(0, \infty ; H_{0}^{1}(\Omega))$ and $u_{m}(t) \in W$ for $0 \leqslant t<\infty$
satisfying
\begin{equation}\label{e32}
 (u_{m,t}, v)+(\nabla u_{m,t}, \nabla v)=2(H(u)u_{m,x}\wedge u_{m,y}, v),\ \ \mbox{for\ all}\  v \in H_{0}^{1}(\Omega), t>0.
\end{equation}
 \begin{equation}\label{e33}
\int_{0}^{t} \int_{\Omega}\left|u_{m,\tau}\right|^{2}+E\left(u_m\left(t\right)\right)
=E(u_{0m})<d, \quad 0 \leqslant t<\infty,
\end{equation}
which implies that
\begin{equation}
E\left(u_{m}\right)=\frac{1}{6}\left\|u_{m}\right\|^{2}+\frac{1}{3} D\left(u_{m}\right).
\end{equation}
So, one has
\begin{equation}
\int_{0}^{T}\left|u_{m \tau}\right|_{2}^{2} d \tau+\frac{1}{6}\left\|u_{m}\right\|_{H_{0}^{1}}^{2}<d, \quad 0 \leqslant t<\infty.
\end{equation}
The remainder of the proof is similar to  the proof of Theorem \ref{t21}.
\end{proof}

\begin{theorem}\label{t32}
(Blow-up) Assume
that $u_{0} \in H_{0}^{1}(\Omega)$, $E\left(u_{0}\right)=d$ and $I\left(u_{0}\right) > 0,$
 Then the existence time of weak solution for system  (\ref{e1}) is finite.
\end{theorem}

\begin{proof}[\textbf{Proof of theorem \ref{t32}}]
Let $u(t)$ be any weak solution of system  \eqref{e1} with $E\left(u_{0}\right)=d$ and $DI\left(u_{0}\right)<0,$ $ T$ be the
existence time of $u(t)$ . We next  prove $T<\infty.$
We argue by contradiction. Suppose that there would exist
a global weak solution $u(t)$. Set
\begin{equation}
f(t)=\int_{0}^{t} \int_{\Omega}|u|^{2}, t>0.
\end{equation}
Multiplying \eqref{e1} by $u$ and integrating over $\Omega \times(0, t),$ we get
\begin{equation}
\int_{\Omega}|u(t)|^{2}-\int_{\Omega}\left|u_{0}\right|^{2}
=-2 \int_{0}^{t} \int_{\Omega}\left(|\nabla u|^{2}+2 H(u) u \cdot u_{x} \wedge u_{y}\right).
\end{equation}
According to the definition of $f(t),$ we have $f^{\prime}(t)=\int_{\Omega}|u(t)|^{2}$ and hence
\begin{equation}
f^{\prime}(t)= \int_{\Omega}|u|^{2}=\int_{\Omega}\left|u_{0}\right|^{2}-2 \int_{0}^{t} \int_{\Omega}\left(|\nabla u|^{2}
+2 H(u) u \cdot u_{x} \wedge u_{y}\right),
\end{equation}
and
\begin{equation}\label{e317}
f^{\prime \prime}(t)=-2 \int_{\Omega}\left(|\nabla u|^{2}+2 H(u) u \cdot u_{x} \wedge u_{y}\right)=-2D(u).
\end{equation}
Now using  \eqref{e12},\eqref{e317}  and
\begin{align}\label{e316}
E(u)&=\frac{1}{2} \int_{\Omega}|\nabla u|^{2}+\frac{2}{3} \int_{\Omega} H(u) u \cdot u_{x} \wedge u_{y}\nonumber\\
&=\frac{1}{6}\|u\|^2+\frac{1}{3}D(u),
\end{align}
we can obtain

\begin{align}
f^{\prime \prime}(t) &= 6 \int_{0}^{t} \int_{\Omega}u^2_{\tau}d\tau
 +f^{\prime}(t)-6 d\nonumber\\
 &=6 \int_{0}^{t} \int_{\Omega}u^2_{\tau}d\tau
 + \int_{\Omega}|u|^{2}-6 d.
\end{align}
Note that

\begin{align}
f(t)f^{\prime \prime}(t) &= f(t)\left[6 \int_{0}^{t} \int_{\Omega}u^2_{\tau}d\tau
 +f^{\prime}(t)-6 E\left(u_{0}\right)\right]\nonumber\\
 &=6\int_{0}^{t} \int_{\Omega}|u|^{2} \int_{0}^{t} \int_{\Omega}u^2_{\tau}d\tau
 +f(t)f^{\prime}(t)-6 d\int_{0}^{t} \int_{\Omega}|u|^{2}.
\end{align}
Hence, we have

\begin{align}\label{3e1}
f(t)f^{\prime \prime}(t) -\frac{3}{2}(f^{\prime}(t))^2&
 \geq 6\int_{0}^{t} \int_{\Omega}|u|^{2} \int_{0}^{t} \int_{\Omega}u^2_{\tau}d\tau-6\left(\int_{0}^{t}\int_{\Omega}u_{\tau}\cdot u \mathrm{d} \tau\right)^{2}\nonumber\\
 &\ \ +f(t)f^{\prime}(t) -3f^{\prime}(t)\int_{\Omega}u_0^2- 6d\int_{0}^{t} \int_{\Omega}|u|^{2}.
\end{align}
Hence, according to \eqref{3e1} and the Schwartz inequality, we obtain
\begin{equation}\label{3e2}
\begin{aligned}
f(t)f^{\prime \prime}(t) -\frac{3}{2}(f^{\prime}(t))^2&
 \geq f(t)f^{\prime}(t) -3f^{\prime}(t)\int_{\Omega}u_0^2- 6 df(t).\\
=&\left(\frac{1}{2} f(t)-3\int_{\Omega}|u_0|^{2}\right) f^{\prime}(t) \\
 &\ \ \ +\left(\frac{1}{2} f^{\prime}(t)-2(p+1) d\right) f(t)>0.
\end{aligned}
\end{equation}
On the other hand, from $E\left(u_{0}\right)=d>0, D\left(u_{0}\right)<0$ and the continuity of $E(u)$ and $D(u)$ with
respect to $t,$ it follows that there exists a sufficiently small $t_{1}>0$ such that $E\left(u\left(t_{1}\right)\right)>0$ and
$D(u)<0$ for $0 \leqslant t \leqslant t_{1} .$
Hence $\left(u_{t}, u\right)=-D(u)>0,|u_{t}|_2>0$  for
$0 \leqslant t \leqslant t_{1}.$  So, using  the continuity
of $\int_{0}^{t}\left|u_{\tau}\right|_{2}^{2} d \tau$, we can choose a $t_{1}$ such
that
\begin{equation}
0<d_{1}=d-\int_{0}^{t_{1}}\left|u_{\tau}\right|_{2}^{2} d \tau<d.
\end{equation}
And by (\ref{e15}), we get
\begin{equation}
0<E\left(u\left(t_{1}\right)\right) = d-\int_{0}^{t_{1}}\left|u_{\tau}\right|_{2}^{2} d \tau=d_{1}<d.
\end{equation}
 So we can choose  $t=t_{1}$ as the initial time, then we obtain  $u(t) \in V_{\delta}$
 for $\delta \in\left(\delta_{1}, \delta_{2}\right), t_{1} \leqslant t<\infty,$ where
$\left(\delta_{1}, \delta_{2}\right)$ is the maximal interval including $\delta=1$
such that $d(\delta)>d_{1}$ for $\delta \in\left(\delta_{1}, \delta_{2}\right) .$ Thus  we
get $D_{\delta}(u)<0$ and $\|u\|>r(\delta)$
for $\delta \in\left(1, \delta_{2}\right), t_{1} \leqslant t<\infty,$
and $D_{\delta_{2}}(u) \leqslant 0,\|u\|\geqslant r\left(\delta_{2}\right)$ for
$t_{1} \leqslant t<\infty .$ Thus  \eqref{e317} implies that
\begin{equation}
\begin{aligned}
 f^{\prime \prime}(t) &=-2 D(u)=2\left(\delta_{2}-1\right)\|\nabla u\|_{2}^{2}-2 D_{\delta_{2}}(u), \\
 & \geqslant 2\left(\delta_{2}-1\right)|\nabla u|_{2}=2\left(\delta_{2}-1\right)\|u\|^{2} \geqslant 2\left(\delta_{2}-1\right) r^{2}\left(\delta_{2}\right), \quad t \geqslant t_1, \\
  f^{\prime}(t) & \geqslant 2\left(\delta_{2}-1\right) r^{2}\left(\delta_{2}\right)(t-t_1)+f^{\prime}(t_1) \geqslant 2\left(\delta_{2}-1\right) r^{2}\left(\delta_{2}\right) (t-t_1), \quad t \geqslant 0,
 \\ f(t) &\geqslant \left(\delta_{2}-1\right) r^{2}\left(\delta_{2}\right) (t-t_1)^{2}+M(t_1)>\left(\delta_{2}-1\right) r^{2}\left(\delta_{2}\right) (t-t_1)^{2}, \quad t \geqslant t_1.
 \end{aligned}
\end{equation}
Therefore, for sufficiently large $t$, we infer
\begin{equation}
\frac{1}{2} f(t)>3\int_{\Omega}|u_0|^{2},\ \  \frac{1}{2} f^{\prime}(t) >6 d.
\end{equation}
Then,  (\ref{3e2}) implies that
\begin{equation}\label{ee2}
\begin{aligned}
f(t)f^{\prime \prime}(t) -\frac{3}{2}(f^{\prime}(t))^2&
 \geq f(t)f^{\prime}(t) -3f^{\prime}(t)\int_{\Omega}u_0^2- 6 E(u_{0})f(t).\nonumber\\
=&\left(\frac{1}{2} f(t)-3\int_{\Omega}|u_0|^{2}\right) f^{\prime}(t) \nonumber\\
 &\ \ \ +\left(\frac{1}{2} f^{\prime}(t)-2(p+1) E\left(u_{0}\right)\right) f(t)>0
\end{aligned}
\end{equation}
The remainder of the proof is the same as that in \cite{r5}.

\end{proof}

\begin{theorem}\label{t33}
Assume that
$u_{0} \in H_{0}^{1}(\Omega)$, $E\left(u_{0}\right)=d$ and $D\left(u_{0}\right)>0,$  $\delta_{1}<\delta_{2}$ are the two roots of equation
$d(\delta)=E(u_0).$  Then, for the global weak solution u of system  \eqref{e1}, it holds
\begin{equation}
|u|_2^2 \leqslant\left|u_{0}\right|_2^{2} e^{-2\left(1-\delta_{1}\right) t}, \quad 0 \leqslant t<\infty.
\end{equation}
\end{theorem}
\begin{proof}[\textbf{Proof of theorem \ref{t33}}]
We first know that system  \eqref{e1} has a global weak solution from Theorem \ref{t32}. Futhermore,
Using  Theorem \ref{t22}, Theorem \ref{t32} and (\ref{e13}), if $u(t)$ is a global weak solution
of system  (\ref{e1}) with $E(u_0) = d$, $D(u_0)> 0$, then must have $D (u)\geq  0$ for $0 \leq t <+\infty$.
Next,  we  distinguish two case:

(1) Suppose that $D(u)>0$ for $0 \leqslant t<\infty$. Multiplying \eqref{e1} by $v$,
 $v \in L^{\infty}\left(0, \infty; H_{0}^{1}(\Omega)\right)$, we have

\begin{equation}\label{e3212}
 (u_t, v)+(\nabla u, \nabla v)=2H(u_x\wedge u_y, v).
\end{equation}
Letting $v = u$, (\ref{e3212}) implies that
\begin{equation}\label{e3213}
\frac{1}{2} \frac{d}{d t}\left|u\right|_2^{2}=-D(u)<0, \quad 0 \leqslant t<\infty.
\end{equation}
Since $\left|u_{t}\right|_{2}>0,$
we have  that $\int_{0}^{t}\left|u_{\tau}\right|^{2} d \tau$ is increasing for $0 \leqslant t<\infty .$ By choosing  any $t_{1}>0$
and letting
\begin{equation} \label{e3214}
d_{1}=d-\int_{0}^{t_{1}}|u_{\tau}|_{2}^{2} d \tau
\end{equation}

From \eqref{e13}, if follows that  $0<E\left(u\right)\leq d_1<d$,
and  $u(t) \in W_{\delta}$ for $\delta_{1}<\delta<\delta_{2}$ and
 $0 \leqslant$ $t<\infty,$ where $\delta_{1}<\delta_{2}$ are the two roots of equation
  $d(\delta)=E\left(u_{0}\right).$ Hence, we obtain  $D_{\delta_1}(u) \geqslant 0$ for
$\delta_{1}<\delta<\delta_{2}$ and $D_{\delta_{1}}(u) \geqslant 0$ for $t_1 \leqslant t<\infty .$ So, \eqref{e3213} gives

\begin{equation}
\frac{1}{2} \frac{d}{d t}\left|u\right|_2^{2}+(1-\delta_1)\left|u\right|_2^{2}+D_{\delta}(u)=0, \quad t_1\leqslant t<\infty.
\end{equation}
Now  \eqref{e3213}  implies that
\begin{equation}
\frac{1}{2} \frac{d}{d t}\left|u\right|_2^{2}+(1-\delta_1)\left|u\right|_2^{2}\leq 0, \quad t_1 \leqslant t<\infty.
\end{equation}
and
\begin{equation}
|u|_2^{2} \leqslant\left|u_{0}\right|_2^{2}-2\left(1-\delta_{1}\right) \int_{0}^{t}|u(\tau)|^{2} d \tau, \quad t_1 \leqslant t<\infty.
\end{equation}
By Gronwall’s inequality, we have

\begin{equation}
|u|_{2}^{2} \leqslant\left|u_{0}\right|_{2}^{2} e^{-2\left(1-\delta_{1}\right) t}, \quad t_1 \leqslant t<\infty.
\end{equation}
\end{proof}

(2) Suppose that there exists a $t_{1}>0$ such that $D\left(u\left(t_{1}\right)\right)=0$ and $D(u)>0$ for $0 \leqslant t<t_{1}$.
Then, $\left|u_{t}\right|_{2}>0$ and $\int_{0}^{t}\left|u_{\tau}\right|_{2}^{2} d \tau$
is increasing for $0 \leqslant t<t_{1}.$  By \eqref{e3214} we have
\begin{equation} \label{e3215}
E(u(t_1))=d-\int_{0}^{t_{1}}|u_{\tau}|_{2}^{2} d \tau<d,
\end{equation}
and $\left\|u\left(t_{1}\right)\right\|=0$. Then, we have that  $u(t) \equiv 0$ for $t_{1} \leqslant t<\infty$.

Hence, the proof is complete.

\section{High initial energy $E(u_0)>d$.}
In this section, we investigate the conditions to ensure the existence of global solutions or blow-up solutions
to system (\ref{e1}) with $E(u_0)>d$.
\begin{lemma}\label{l59}
For any $\alpha>d, \lambda_{\alpha}$ and $\Lambda_{\alpha}$ defined in \eqref{ee123} satisfy
\begin{equation}
0<\lambda_{\alpha} \leq \Lambda_{\alpha}<+\infty.
\end{equation}
\end{lemma}
\begin{proof}
(1) By  H\"{o}lder inequality,   fundamental inequality and $u\in \mathcal{N}$, we have
\begin{equation}
\begin{aligned}
\frac{1}{2} \int_{\Omega}|\nabla u|^{2}&=\left|\int_{\Omega} H(u) u \cdot u_{x} \wedge u_{y}\right| \\
 &\leq H\left(\int_{\Omega} u^2\right)^{\frac{1}{2}}
\left(\int_{\Omega}  |  u_{x} \wedge u_{y}|^2\right)^{\frac{1}{2}}\\
&\leq CH\left(\int_{\Omega} u^2\right)^{\frac{1}{2}}
\left(\int_{\Omega}|\nabla u|^{2}\right)^{\frac{1}{2}}.
\end{aligned}
\end{equation}
Then from Lemma \ref{l2801} (1), we have $\lambda_{\alpha}>0$.

(2) Using the isoperimetric inequality and $u\in \mathcal{N}$, we have
 \begin{equation}
 \frac{1}{2|H|}\int_{\Omega}|\nabla u|^{2}
 =\left|\int_{\Omega} u \cdot u_{x} \wedge u_{y}\right|\leq \frac{1}{4\sqrt{2\pi}}\left(\int_{\Omega}|\nabla u|^{2}\right)^{\frac{3}{2}}.
\end{equation}
So we have $\|u\|\leq \frac{|H|}{2\sqrt{2\pi}},$  which  leads to the conclusion.
\end{proof}

\begin{theorem}\label{t51}
Suppose  that $E\left(u_{0}\right)>d$, then we have
\begin{itemize}\addtolength{\itemsep}{-1.5 em} \setlength{\itemsep}{-5pt}
\item[\emph{(1)}]  If $u_{0} \in \mathcal{N}_{+}$ and $\left|u_{0}\right|_{2} \leq \lambda_{E\left(u_{0}\right)},$ then $u_{0} \in \mathcal{G}_{0}$,
\item[\emph{(2)}] If $u_{0} \in \mathcal{N}_{-}$ and  $\left|u_{0}\right|_{2} \geq \Lambda_{E\left(u_{0}\right)},$ then $u_{0} \in \mathcal{B}$.
\end{itemize}
\end{theorem}
\begin{proof}
 The maximal existence time    of the solutions to system  (\ref{e1}) with initial value
$u_0$ is  denoted by $T_0$.
 If the solution is global, i.e. $T(u_0) =+\infty$,   the  limit set of $u_0$ is denoted by $\omega_0$.

 (1) Suppose  that $u_{0} \in \mathcal{N}_{+}$ with $\left|u_{0}\right|_{2} \leq \lambda_{E\left(u_{0}\right)}.$
 We firstly prove that $u(t) \in \mathcal{N}_{+}$ for all $t \in\left[0, T\left(u_{0}\right)\right) .$
Assume, on the contrary, that there exists a $t_{0} \in\left(0, T\left(u_{0}\right)\right)$ such that $u(t) \in \mathcal{N}_{+}$ for $0 \leq t<t_{0}$
 and $u\left(t_{0}\right) \in \mathcal{N} .$ It follows from
$D(u(t))=-\int_{\Omega} u_{t}(x, t) u(x, t) \mathrm{d} x$ that $u_{t}(x, t) \neq 0$ for $(x, t) \in \Omega \times\left(0, t_{0}\right) .$
Recording to\eqref{e12} we then have
$E\left(u\left(t_{0}\right)\right)<E\left(u_{0}\right),$ which implies that $u\left(t_{0}\right) \in E^{E\left(u_{0}\right)} .$
 Therefore, $u\left(t_{0}\right) \in \mathcal{N}^{E\left(u_{0}\right)}.$    Recalling  the definition
of $\lambda_{E\left(u_{0}\right)},$ we get
  \begin{equation}\label{e51}
\left|u\left(t_{0}\right)\right|_{2} \geq \lambda_{E}\left(u_{0}\right).
\end{equation}
Since  $D(u(t))>0$ for $t \in\left[0, t_{0}\right),$  we obtain from \eqref{e213} that
\begin{equation}
\left|u\left(t_{0}\right)\right|_{2}<\left|u_{0}\right|_{2} \leq \lambda_{E\left(u_{0}\right)}.
\end{equation}
which contradicts \eqref{e51}.
Hence,  $u(t) \in \mathcal{N}_{+}$ which shows  that $u(t) \in E^{E\left(u_{0}\right)}$ for all
$t \in\left[0, T\left(u_{0}\right)\right) .$  Now  Lemma \ref{l28} (2) implies that the orbit $\{u(t)\}$
remains bounded in $H_{0}^{1}(\Omega)$ for $t \in\left[0, T\left(u_{0}\right)\right)$
so that $T\left(u_{0}\right)=\infty .$
 Assume that  $\omega$ is an arbitrary element in $\omega\left(u_{0}\right).$ Then by (\ref{e12}) and \eqref{e213} we obtain
\begin{equation}
|\omega|_{2}>\Lambda_{E\left(u_{0}\right)}, \quad E(\omega)<E\left(u_{0}\right),
\end{equation}
which, according to  the definition of $\lambda_{E\left(u_{0}\right)}$ again, implies that  $\omega\left(u_{0}\right) \cap N=\emptyset .$ So, $\omega\left(u_{0}\right)=\{0\},$ i.e. $u_{0} \in \mathcal{G}_{0}$

(2)  Suppose  that $u_{0} \in \mathcal{N}_{-}$ with $\left|u_{0}\right|_{2} \geq \Lambda_{E\left(u_{0}\right)}.$
 We now  prove that $u(t) \in \mathcal{N}_{-}$ for all $t \in\left[0, T\left(u_{0}\right)\right) .$
 Assume, on the contrary, that
there exists a $t^{0} \in\left(0, T\left(u_{0}\right)\right)$ such that $u(t) \in \mathcal{N}_{-}$
 for $0 \leq t<t^{0}$ and $u(t^{0}) \in \mathcal{N} .$
Similarly to case (1), one has
$E(u(t^{0}))<E\left(u_{0}\right),$
which implies that $u(t^{0}) \in E^{E\left(u_{0}\right)} .$
Therefore, $u(t^{0}) \in \mathcal{N}^{E\left(u_{0}\right)} .$
Recalling the definition of $\Lambda_{E\left(u_{0}\right)},$ we infer
\begin{equation}\label{e52}
|u(t^{0})|_{2} \leq \Lambda_{E\left(u_{0}\right)}
\end{equation}
On the other hand, from \eqref{e213} and the fact that $D(u(t))<0$ for $t \in[0, t^{0}),$ we obtain
\begin{equation}
|u(t^{0})|_{2}>|u_{0}|_{2} \geq \Lambda_{E\left(u_{0}\right)},
\end{equation}
which contradicts \eqref{e52}.

Assume that $T\left(u_{0}\right)=\infty .$ Then for each $\omega \in \omega\left(u_{0}\right),$
it follows from by (\ref{e12}) and \eqref{e213} that
\begin{equation}
\|\omega\|_{2}>\Lambda_{E\left(u_{0}\right)}, \quad E(\omega)<E\left(u_{0}\right).
\end{equation}
Noting the definition of $\Lambda_{E\left(u_{0}\right)}$ again, we have  $\omega\left(u_{0}\right) \cap N=\emptyset .$
 Hence, it is holded that $\omega\left(u_{0}\right)=\{0\}$,
which contradicts Lemma \ref{l28} (1). Therefore, $T\left(u_{0}\right)<\infty$ and we can complete this proof.

\end{proof}

\begin{theorem}\label{t52}
Assume that $u_{0} \in H_{0}^{1}(\Omega)$ satisfies
\begin{equation}\label{e54}
E(u_0)\leq|u_0|_2<-\frac{1}{3}\int_{\Omega} H(u_0) u_0 \cdot u_{0,x} \wedge u_{0,y},
\end{equation}
Then$, u_{0} \in \mathcal{N}_{-} \cap \mathcal{B}$.
\end{theorem}
\begin{proof}
Firstly, we observe
\begin{equation}
\begin{aligned}
E(u_0)&=\frac{1}{2} \int_{\Omega}|\nabla u_0|^{2}+\frac{2}{3} \int_{\Omega} H(u_0) u_0 \cdot u_{0,x} \wedge u_{0,y} \\
&= \frac{1}{2}D(u_0)-\frac{1}{3} \int_{\Omega} H(u_0) u_0 \cdot u_{0,x} \wedge u_{0,y}
 \end{aligned}
\end{equation}
Thus, we have
\begin{equation}\label{e53}
 E(u_0)+\frac{1}{3} \int_{\Omega} H(u_0) u_0 \cdot u_{0,x} \wedge u_{0,y}=\frac{1}{2}D(u_0)<0,
\end{equation}
which shows that $u_0\in \mathcal{N}_{-}$.
Then for any $u \in \mathcal{N}_{E\left(u_{0}\right)},$   one has
$$-2\int_{\Omega} H(u) u \cdot u_{x} \wedge u_{y}=\|u\|^2\leq 6E(u_0).$$
Taking supremum over $\mathcal{N}_{E\left(u_{0}\right)}$ and \eqref{e54}, by Theorem \ref{t51}  we can deduce
$$\left|u_{0}\right|_{2} \geq \Lambda_{E\left(u_{0}\right)}.$$
 Thus, $u_{0} \in \mathcal{N}_{-} \cap \mathcal{B}$.
\end{proof}





\section*{Acknowledgements.}
This  work  is supported by   NSFC (No. 11801017)  and Beijing Municipal Natural Science Foundation  (No. 1172005).

\end{document}